\documentclass{amsart}

\newtheorem{theorem}{Theorem}
\newcommand{\bt}{\begin{theorem}}
\newcommand{\et}{\end{theorem}}
\newtheorem{lemma}{Lemma}
\newcommand{\bl}{\begin{lemma}}
\newcommand{\el}{\end{lemma}}
\newtheorem{corollary}{Corollary}
\newcommand{\bc}{\begin{corollary}}
\newcommand{\ec}{\end{corollary}}
\newtheorem{problem}{Problem}
\newcommand{\bprob}{\begin{problem}}
\newcommand{\eprob}{\end{problem}}
\newcommand{\beq}{\begin{equation}}
\newcommand{\eeq}{\end{equation}}
\newcommand{\benum}{\begin{enumerate}}
\newcommand{\eenum}{\end{enumerate}}

\newcommand{\N}{\ensuremath{ \mathbf N }}
\newcommand{\Z}{\ensuremath{\mathbf Z}}

\newcommand{\R}{\ensuremath{\mathbf R}}

\newcommand{\mbA}{\ensuremath{ \mathbf A}}

\newcommand{\mbF}{\ensuremath{ \mathbf F}}
\newcommand{\PP}{\ensuremath{ \mathbf P}}

\newcommand{\mcc}{\ensuremath{ \mathcal C}}
\newcommand{\mcd}{\ensuremath{ \mathcal D}}

\newcommand{\mcf}{\ensuremath{ \mathcal F}}
\newcommand{\mcg}{\ensuremath{ \mathcal G}}

\newcommand{\mcw}{\ensuremath{ \mathcal W}}

\newcommand{\mba}{\ensuremath{ \mathbf a}}

\newcommand{\mbh}{\ensuremath{ \mathbf h}}

\DeclareMathOperator{\card}{\text{card}}

\newcommand{\bmat}{\left(\begin{matrix}}
\newcommand{\emat}{\end{matrix}\right)}

\DeclareMathOperator{\qqand}{\qquad\text{and}\qquad}

\title[Extremal problems]{Extremal problems and the combinatorics of sumsets}

\author{Melvyn B. Nathanson}
\address{Lehman College (CUNY), Bronx, New York 10468}
\email{melvyn.nathanson@lehman.cuny.edu}

\date{\today}



\thanks{Supported in part by grant \# 66197-00 54 from the PSC-CUNY Research Award Program.}

\dedicatory{To Helmut Maier on his seventieth birthday}

\begin{document}
\maketitle

\begin{abstract} 
This is a survey of old and new problems and results in additive number theory. 
\end{abstract}

\section{Sumsets and asymptotic bases}

Here are the best known theorems in additive number theory.
Every nonnegative integer is the sum of four squares (Lagrange).   
Every nonnegative integer is the sum of nine nonnegative cubes (Wieferich), 
every sufficiently large integer is the sum of seven nonnegative cubes (Linnik), 
and  almost every  nonnegative integer is the sum 
of four nonnegative cubes (Davenport).  
For all $k \geq 2$ there is a number $g(k)$ such that 
every nonnegative integer is the sum of $g(k)$ $k$th powers (Hilbert). 
Every sufficiently large odd integer is the sum of three primes (Vinogradov).  
These results  
are statements about sums of a fixed number of elements chosen from 
an ``interesting'' set of nonnegative integers.  
For a survey of these results, see Nathanson~\cite{nath1996aa}.  

The classical theorems motivate the following fundamental definitions.
Let $A, A_1,\ldots, A_h$ be finite or infinite subsets of the set $\N_0 = \{0,1,2,3,\ldots\}$ 
of nonnegative integers.  Let  $\N = \{1,2,3,\ldots\}$ be the set of positive integers.
We define the \emph{sumset} 
\[
A_1 + \cdots + A_h = \{a_1 + \cdots + a_h: a_i \in A_i \text{ for all } i  = 1,2,\ldots, h \} 
\]
and the \emph{$h$-fold sumset} 
\[
hA = \underbrace{A + \cdots + A}_{\text{$h$ summands}} = \{a_1 + \cdots + a_h: 
a_i \in A \text{ for all } i =1,2,\ldots, h\}.
\]
We also consider sums of subsets of the set $\Z^n$ of  lattice points, 
of the the set $\Z/m\Z$  congruence classes modulo $m$, 
and, more generally, of any additive abelian semigroup.
  
The set $A$ of nonnegative integers is a \emph{basis of order $h$} 
if every nonnegative integer can be represented 
as the sum of $h$ not necessarily distinct elements of $A$, that is, if $hA = \N_0$.  
The set $A$ is an \emph{asymptotic basis of order $h$} if $hA$ contains all sufficiently large integers, or, 
equivalently, if $hA \sim \N_0$.  
(We write $X \sim Y$ if the symmetric difference $(X\setminus Y) \cup (Y \setminus X)$ is finite.)    
Thus, Lagrange's theorem states that the set $A^{(2)} = \{k^2: k \in \N_0\}$ 
is a basis of order 4,  Wieferich's theorem states that the set $A^{(3)} = \{k^3: k \in \N_0\}$ 
is a basis of order 9, and Linnik's theorem states that $A^{(3)}$ is an asymptotic basis of order 7. 

Let $A$ be a set of nonnegative integers.  The {$h$-fold representation function} 
$r_{h,A}(n)$ counts the number of representations of $n$ as a sum of $h$ elements of $A$, 
that is, the number of $h$-tuples 
\[
(a_1,\ldots, a_h) \in A^h 
\] 
such that 
\[
n = a_1 + a_2 + \cdots + a_h 
\]
and 
\[
a_1 \leq a_2 \leq \cdots \leq a_h.
\]
A central unsolved problem in additive number theory is the following conjecture 
of Erd\" os and  Tur\" an~\cite{erdo-tura41}.

\bprob 
Let $A$ be a set of nonnegative integers.  
If $A$ is an asymptotic basis of order $h$, 
then the representation function $r_{h,A}(n)$ is unbounded.  Equivalently, 
\[
\liminf_{n\rightarrow \infty} r_{h,A}(n) > 0 
\]
implies 
\[
\limsup_{n\rightarrow \infty} r_{h,A}(n) = \infty. 
\]
\eprob

One approach to the Erd\H os-Tur\' an conjecture is to look for 
extremal asymptotic bases, that is, bases  
that are, in some sense, ``small''  and might have small representation functions. 
We consider two kinds of extremal bases: Thin bases and minimal bases.  
This paper discusses old and new results and old and new problems 
related to these extremal constructions in combinatorial and additive number theory.

The Erd\" os-Tur\" an conjecture was originally stated for bases of order 2, 
but extends naturally to bases of all orders $h \geq 2$. 
Many theorems in additive number theory about 2-fold sumsets  
do not extend or have not (yet) been extended to $h$-fold sumsets for $h \geq 3$.  
One difficulty that  often arises in combinatorial and probabilistic arguments is 
the lack of independence of representations.  
Two representations 
\[
n = a_1 +\cdots + a_h \qqand n = a'_1 +\cdots + a'_h
\]
are \emph{independent} if 
\[
\{a_1,\ldots, a_h\} \cap \{a'_1,\ldots, a'_h\} = \emptyset.
\]
All representations with $h=2$ summands are independent:  
If $n = a_1+a_2=a'_1+a'_2$, 
then either  $\{a_1, a_2\} = \{ a'_1, a'_2\}$ or $\{a_1, a_2\} \cap \{ a'_1, a'_2\} = \emptyset$.  
This is not true for representations as sums of more than two elements; 
for example, $8 = 1 + 2 + 5 = 1 + 3 + 4$.  
This phenomenon, known for a long time, has recently been rebranded as ``additive energy.''

In the study of $h$-fold sumsets of a set $A$, we define the 
\emph{solution set} $S_{h,A}(n)$:
\[
S_{h,A}(n) = \{a_1 \in A: \text{there exist $a_2,\ldots,a_h \in A$ such that 
$n = a_1+a_2+\cdots + a_h$} \}.
\]

Let $\Omega$ be the space of all infinite sets of positive integers, and let 
$p_n \in [0,1]$ for all $n \in \N$.  Let 
\[
E_n = \{A\in \Omega: n \in A\}.
\]
Erd\H os and R\' enyi~\cite{erdo-reny60} constructed a probability measure $\mu$ 
on $\Omega$ such that $\mu(E_n) = p_n$ 
and the events $E_1, E_2, E_2,\ldots$ are independent.
We have the following result about solution sets.

\bt[Erd\H os-Nathanson~\cite{nath1986-54}]  \label{nsf:theorem:nath1986-54}
For $\alpha, \beta, \gamma \in \R$ with $\alpha > 0$ 
and $1/3 < \gamma \leq 1/2$, let 
\[
p_n \leq \frac{\alpha \log^{\beta}(n+1)}{n^{\gamma}} 
\]
for all $n \in \N$.  
For almost all $A \in \Omega$ and for all but finitely many 
pairs $(m,n)$ with $m \neq n$, 
\[
|S_{2,A}(m) \cap S_{2,A}(n)| \leq \frac{2}{3\gamma-1}.
\]
\et

Choosing $\gamma = 1/2$, we obtain sets $A$ such that 
$|S_{2,A}(m) \cap S_{2,A}(n)| \leq 4$ 
for all but finitely many pairs $(m,n)$ of integers $m \neq n$.

\bprob
Extend Theorem~\ref{nsf:theorem:nath1986-54} to solution sets 
$S_{h,A}(n)$ for $h \geq 3$.
\eprob

\section{Thin bases}

Our first measure of smallness is ``slow growth.''
Let $A$ be a set of nonnegative integers.  
The \emph{counting function} $A(x)$  counts the number of 
elements of $A$ that are no greater than $x$,
that is,
\[
A(x) = \sum_{\substack{a\in A\\ a \leq x}} 1. 
\]
If $A$ is an asymptotic basis of order $h$, then $r_{h,A}(n) \geq 1$ for all $n \geq n_0$.  
If $n_0 \leq n \leq x$ and $n = a_1+\cdots + a_h \in hA$, then $a_i \leq x$ for all 
$i \in \{1,\ldots, h\}$.  The number of $h$-tuples $(a_1,\ldots, a_h) \in A^h$ 
with $a_1 \leq \cdots \leq a_h \leq x$ is $\binom{A(x)+h-1}{A(x)-1} = \binom{A(x)+h-1}{h}$.   
It follows that 
\[
x-n_0 \leq \binom{A(x)+h-1}{h} 
\]
and so 
\[
A(x) \gg x^{1/h}.
\]
An asymptotic basis $A$ of order $h$ is \emph{thin} if 
\[
A(x) \ll x^{1/h}.
\]
There are simple constructions of thin bases due to Raikov~\cite{raik37}, Jia~\cite{jia90,jia92,jia96}, and Nathanson~\cite{nath1974-11,nath1988-66,nath2012-145}.  
Shatrovskii~\cite{chat40,nath2012-145} constructed beautiful examples of thin bases
and, for every $h \geq 2$, Cassels~\cite{cass57,nath2010-138} constructed a thin 
asymptotic basis 
$A = \{a_k:k=1,2,\ldots\}$ 
such that $a_k = ck^h + O\left(k^{h-1}\right)$ for some $c > 0$.

An asymptotic basis of order $h$ does not necessarily contain a subset 
that is a thin basis of order $h$, but it is natural to ask if the classical bases 
in additive number theory contain thin bases.

The set $A^{(2)} = \{k^2:k \in \N_0\}$ is a basis of order 4 
with counting function $A(x) \sim x^{1/2}$.  A basis of order 4 must have counting 
function  $\gg x^{1/4}$ and so the set of all squares is  ``thick''.   
It is an unsolved problem 
to determine if there is a set of squares 
that is a thin basis of order 4.  Probability arguments do show that there 
exist sets  $W$ of nonnegative integers of density 0 such that 
the set of squares $\{w^2: w \in W\}$ is a basis of order 4. 
No explicit example is known.  For work on this problem, see 
Choi-Erd\H os-Nathanson~\cite{nath1980-44}, 
Erd\H os-Nathanson~\cite{nath1981-46}, Wirsing~\cite{wirs86}, 
and Z\" ollner~\cite{zoll84,zoll85}. 

\bprob
Construct a set $W$ of nonnegative integers 
of density 0 such that the set $\{w^2:w\in W\}$ is a basis of order 4. 
\eprob

There is the analogous problem for sums of $k$th powers (Waring's problem).  
Hardy and Littlewood proved  
that for every $k \geq 2$ there is an integer 
$s_0(k)$ such that, if $s \geq s_0(k)$, then a positive integer has $cn^{(s/k)-1}$ 
representations as a sum of $k$th powers 
(Davenport~\cite{dave62}, Vaughan~\cite{vaug97}). This is used to prove the following. 

\bt[Nathanson~\cite{nath1981-47}]
For any $s > s_0(k)$ and $0 < \varepsilon < 1/s$, there is a probability 
measure on the set $\Omega$ of all infinite sets of positive integers such that, 
with probability 1, a random set $W \in \Omega$ has the following properties:
\benum
\item[(i)]
$W(x) \sim cx^{1 - (1/s)+\varepsilon}$ for some $c > 0$,
\item[(ii)]
every sufficiently large integer $n$ can be repsented in the form 
\[
n = w_1^k + 2_2^k + \cdots + w_s^k 
\]
with $w_1,w_2,\ldots, w_s \in W$. 
\eenum
\et  

Van Vu~\cite{vu00} has improved this result. 

\bprob
Construct a set $W_k$ of nonnegative integers of density 0 
such that the set $\{w^k:w\in W_k\}$ 
is a basis of order h for some $h$. 
\eprob

\section{Minimal asymptotic bases}
A second measure of the smallness of an asymptotic basis is ``minimality''. 
An asymptotic basis $A$ of order $h$ is \emph{minimal} if no proper subset of $A$ 
is an asymptotic basis of order $h$.  This means that, for every $a \in A$, 
the set $A\setminus \{a\}$ is not an asymptotic basis of order $h$, and so 
$h\left( A\setminus \{a\} \right)$ is a co-infinite subset of $\N_0$.  
Thus, for every $a \in A$, there are infinitely many nonnegative integers $n$,  
all of whose representations as sums of $h$ elements of $A$ 
must include the integer $a$.  

We ask the following questions: 
Do minimal  asymptotic bases exist, and, if so,  does every  asymptotic  basis 
contain   a minimal  asymptotic nonbasis?    
In fact, minimal asymptotic bases of order $h$ do exist for all $h \geq 2$, 
but there are few explicit constructions.  
Nathanson~\cite{nath1974-11,nath1988-66} used a 2-adic construction to produce 
the first examples of minimal asymptotic bases.    
This method was extended to $g$-adically defined sets by  
Chen~\cite{chen12}, 
Chen and Chen~\cite{chen-chen11}, Chen and Tang~\cite{chen-tang18}, 
Jia~\cite{jia96}, Jia and Nathanson~\cite{nath1989-68}, 
Lee~\cite{lee93}, Li and Li~\cite{li-li16}, 
Ling and Tang~\cite{ling-tang15,ling-tang18}, 
Sun~\cite{sun19,sun21}, and 
Sun and Tao~\cite{sun-tao19}.  
In these examples, one has $r_{h,A}(n) = 1$ for infinitely many $n$, 
but also $\limsup_{n\rightarrow \infty} r_{h,A}(n) = \infty$, 
so the known examples do not provide counterexamples to the Erd\H os-Tur\" an conjecture. 

Nathanson~\cite{nath2022-201} recently constructed a new class of minimal asymptotic bases 
of order $h$ for all $h \geq 2$.  
A \emph{\mcg-adic sequence} is a strictly increasing sequence of positive integers 
$\mcg = (g_j)_{j=0}^{\infty}$ 
such that $g_ 0 = 1$ and $g_{j-1}$ divides $g_j$ for all $ j \geq 1$.   
Let $W$ be a nonempty set of nonnegative integers, 
and let $\mcf^*(W)$ be the set of all nonempty finite subsets of $W$.  
Let $\mcg = (g_j)_{j=0}^{\infty}$ be a \mcg-adic sequence with $d_{j+1} = g_{j+1}/g_j$. 
We define the set of positive integers 
\[
A_{\mcg}(W) = \left\{ \sum_{j\in F} x_j g_j : F \in \mcf^*(W) \text{ and }  x_j \in [1,d_{j+1} -1] \right\}.
\]

\bt[Nathanson~\cite{nath2022-201}]                  \label{NSF:theorem:nath2022-201} 
Let $h \geq 2$ and let  $t$ be an integer such that 
\[
t \geq 1 +  \frac{\log h}{\log 2}.
\]
Let $\mcw = (W_i)_{i=0}^{h-1}$ be a partition of $\N_0$  such that, 
for all $i \in [0 ,h-1]$, the set $W_i$ contains infinitely many intervals 
of $t$ consecutive integers.  
Let $\mcg$ be a \mcg-adic sequence. 
The set 
\[
A_{\mcg}(\mcw) = \bigcup_{i=0}^{h-1} A_{\mcg}(W_i)
\]
is a minimal asymptotic basis of order $h$.  
\et

Nathanson~\cite{nath1974-11} proved that not every asymptotic basis contains 
a minimal asymptotic basis.  There is the following remarkable example.

\bt[Erd\H os-Nathanson~\cite{nath1975-13}]
There exists an asymptotic basis $A$ of order 2 such that, if $S$ is a subset of $A$, 
then $A\setminus S$ is an asymptotic basis of order 2 if and only if 
the set $S$ is finite.
\et

\bprob
Let $h \geq 3$.  
Construct an asymptotic basis $A$ of order $h$ such that, if $S$ is a subset of $A$, 
then $A\setminus S$ is an asymptotic basis of order $h$ if and only if 
the set $S$ is finite.
\eprob
 
It is a difficult problem to determine if an asymptotic basis contains 
a minimal asymptotic basis.  Here is a result that uses the representation function. 

\bt[Erd\H os-Nathanson~\cite{nath1979-38}]          \label{NSF:theorem:4/3}
Let $A$ be an asymptotic basis of order $2$.
If 
\beq            \label{NSF:c-log-n}
r_{2,A}(n) > \frac{1}{\log(4/3)} \log n
\eeq
for all sufficiently large $n$, then $A$ contains a minimal asymptotic basis of order 2. 
\et

\bprob
The dependence on the particular constant $1/\log(4/3)$ in Theorem~\ref{NSF:theorem:4/3} 
is not understood.  
It is open problem to determine if the condition 
\[
r_{2,A}(n) > c \log n
\]
for some $c > 0$ and all sufficiently large $n$ implies that 
$A$ contains a minimal asymptotic basis of order 2.  
\eprob

\bprob
Determine if the condition 
\[
\lim_{n\rightarrow \infty} r_{2,A}(n) = \infty 
\]
implies that  $A$ contains a minimal asymptotic basis of order 2.
\eprob

\bprob
Let $h \geq 3$.  
Does there exists a constant $c_h > 0$ 
such that the condition 
\[
r_{h,A}(n) > c_h \log n
\]
for all sufficiently large $n$ implies that 
$A$ contains a minimal asymptotic basis of order $h$?
\eprob

Inequality~\eqref{NSF:c-log-n} appears mysteriously in the following result about partitioning 
an asymptotic basis into a disjoint union of asymptotic bases. 

\bt[Erd\H os-Nathanson~\cite{nath1988-63}]
Let $A$ be an asymptotic basis of order 2 whose representation function 
satisfies~\eqref{NSF:c-log-n}.  There are subsets $A_1$ and $A_2$ of $A$ 
such that  $A = A_1 \cup A_2$ and  $A_1 \cap A_2 = \emptyset$ 
and both $A_1$ and $A_2$ are asymptotic bases 
of order $h$. 
\et

Again, the dependence on the constant $1/\log(4/3)$ is not understood.

\section{Maximal asymptotic nonbases}
Minimal asymptotic bases are extremal objects in additive number theory. \\
Nathanson~\cite{nath1974-11} introduced the dual notion.  
A set of nonnegative integers is an \emph{asymptotic nonbasis of order $h$} 
if it is not an asymptotic basis of order $h$.   
An asymptotic nonbasis $A$ of order $h$ is \emph{maximal} if $A \cup \{a\}$
is an asymptotic basis of order $h$ for every nonnegative integer $a \notin A$.
We ask the following questions: 
Do maximal  asymptotic nonbases exist, and, if they do,  is every  asymptotic nonbasis 
contained in a maximal  asymptotic nonbasis? 

There are  trivial examples of maximal asymptotic nonbases. 
For example, the set of even integers is a maximal asymptotic nonbasis of order 2.  
Other simple examples can be constructed that are unions of congruence classes.
There are also nontrivial examples.  Erd\H os and Nathanson~\cite{nath1975-13} 
constructed a class of maximal  asymptotic nonbases of order 2 that 
contain arbitrarily long gaps and so are 
not unions of congruence classes.

The first example of an asymptotic  nonbasis that is not contained in a maximal 
 asymptotic nonbasis is due to 
Hennefeld~\cite{henn77}, who constructed a set $A$ with the property that a superset $A'$ 
of $A$ is a nonbasis if and only if $A'$ is a co-infinite subset of $\N_0$.  
There exists no maximal co-infinite subset of an infinite set, and so $A$ is not contained 
in a maximal nonbasis. 

Here is a recent construction of a new class of asymptotic nonbases that are not 
contained in maximal asymptotic nonbases. 
The set $Y$  of integers has \emph{infinite gaps}\index{infinite gaps} 
if, for all $C > 0$, 
there exist only finitely many pairs of integers $y,y' \in Y$ 
with $y \neq y'$ and $|y-y'| \leq C$.

\bt[Nathanson~\cite{nath2020-188}]                        \label{NSF:theorem:nath2020-188}
Let $h \geq 2$.    
Let $s$ and $t$ be integers such that $\gcd(h,s-t) = 1$. 
\benum
\item[(a)]
Let $Y$ be an infinite set of integers with infinite gaps, and let
$X = \Z \setminus Y.$
The set 
\[
A_X = \{s\} \cup \{hx+t: x \in X \}
\]
is an asymptotic nonbasis of order $h$ for \Z\ 
that is not a subset of a maximal asymptotic nonbasis of order $h$ for \Z. 

\item[(b)]
Let $Y_0$ be an infinite set of  nonnegative integers with infinite gaps, and let
$X_0 = \N_0 \setminus Y_0$. 
The set  
\[ 
A_{X_0} = \{s\} \cup \{hx+t: x \in X_0 \}
\]
is an asymptotic nonbasis of order $h$ for $\N_0$ 
that is not a subset of a maximal asymptotic nonbasis of order $h$ for ${\N_0}$. 
\eenum
\et

\section{Oscillations}

We can extend the notion of minimal basis and maximal nonbasis by considering 
perturbations of a set $A$ of nonnegative integers by a finite 
set $F$ of integers with more than one element.  
Let $A$ be an asymptotic basis of order 2, and let $F$ be a  subset of $A$.  
The set $A$ is \emph{$r$-minimal} if $A\setminus F$ is an asymptotic  basis 
if $|F| < r$ and an asymptotic  nonbasis if $|F| \geq r$.
An asymptotic basis $A$ is $\aleph_0$-minimal if $A\setminus F$ is a basis 
if $F$ is finite and a nonbasis if  $F$ is infinite.

Let $A$ be an asymptotic nonbasis of order 2, 
and let $G$ be a  subset of $\N_0 \setminus A$.  
The set $A$ is \emph{$s$-maximal} if $A\cup G$ is an asymptotic  nonbasis 
if $|G| < s$ and an asymptotic  basis if $|G| \geq s$.
An asymptotic nonbasis $A$ is $\aleph_0$-maximal if $A \cup G$ is a nonbasis 
if $G$ is finite and a  basis if  $G$ is infinite.

\bt[Erd\H os-Nathanson~\cite{nath1975-17}] 
Let $Q = \{2q_k+1\}_{k=1}^{\infty}$ be a set of odd integers such that
$q_k \geq 5q_{k-1}+3$
for all $k \geq 2$.  Let
\[
A^Q = \bigcup_{k=2}^{\infty} \left( \left[ 2q_{k-1}+2, q_k - q_{k-1} \right] 
\cup  \left[ q_k + 1, q_k + q_{k-1} \right] \right).
\]

\benum
\item
For all $r \geq 1$, the set $A^Q$ is contained in an $r$-minimal basis of order 2.

\item
The set $A^Q$ is contained in an $\aleph_0$-minimal basis of order 2.

\item
For all $s \geq 1$, the set $A^Q$ is contained in an $s$-maximal nonbasis of order 2.

\eenum
\et

\bt[Erd\H os-Nathanson~\cite{nath1975-17}] 
\benum
\item
There does not exist an asymptotic  basis $A = \{a_i:i=1,2,3,\ldots\}$ 
 of order 2 such that
$A \setminus \{a_u:u \in U\}$ is an asymptotic  basis  of order 2 
if $U$ has zero density 
and an asymptotic  nonbasis  of order 2 if $U$ has positive density.

\item
There does not exist an $\aleph_0$-maximal asymptotic nonbasis  of order 2.
\eenum
\et

A set of positive integers is \emph{infinitely oscillating} if it switches from basis to 
nonbasis to basis to nonbasis \ldots of order 2 as integers are alternately adjoined to 
and deleted from the set. 
Equivalently, let $A$ be a set of positive integers and let $F$ and $G$ be finite sets 
such that $F \subseteq A$ and $G \subseteq \N\setminus A$. 
The set $A$ is an \emph{infinitely oscillating basis} if 
$(A\setminus F)\cup G$ is an asymptotic basis of order 2 if and only is $|F| \leq |G|$.  

Let $B$ be a set of positive integers and let $F$ and $G$ be finite sets 
such that $G \subseteq B$ and $F \subseteq \N\setminus B$. 
The set $B$ is an infinitely oscillating nonbasis if 
$(B\cup F)\setminus G$ is an asymptotic nonbasis of order 2 if and only is $|F| \leq |G|$. 

Infinitely oscillating bases exist.  
The following stronger result is strange and fascinating.  

\bt[Erd\H os-Nathanson~\cite{nath1976-22}]  \label{NSF:theorem:nath1976-22}
There exists a partition of the positive integers \N\ into two disjoint sets 
$A$ and $B$ such that $A$ is an infinitely oscillating basis 
and $B$ is an infinitely oscillating nonbasis. 
\et

The proof  is nonconstructive.

\bprob
Construct  a partition of the positive integers \N\ into two disjoint sets 
$A$ and $B$ such that $A$ is an infinitely oscillating asymptotic basis of order 2 
and $B$ is an infinitely oscillating asymptotic nonbasis of order 2. 
\eprob

Theorem~\ref{NSF:theorem:nath1976-22} suggests several new problems. 
If $A$ is an asymptotic basis of order $h$, 
then $A$ is an asymptotic basis of order $h'$ for all $h' \geq h$. 
The \emph{exact order} of a set $A$ of integers is the least integer $h$ such that 
$A$ is an asymptotic basis of order $h$.

\bprob
Let $h_0 < h_1 < h_2 < \cdots$  be a strictly increasing 
finite or infinite sequence 
of nonnegative integers with $h_0 \geq 2$.  
\benum

\item
Does there exists an asymptotic basis $A$ of order $h_0$ 
and an infinite sequence $(a_i)_{i=1}^{\infty}$ of elements of $A$ 
such that the set $A\setminus \{a_1,a_2,\ldots, a_i \}$ has exact order $h_i$ for all $i$?
\item
Does there exist an asymptotic basis $A$ of order $h_0$ 
and an infinite sequence $(a_i)_{i=1}^{\infty}$ of elements of $A$ 
such that the set $A\setminus \{a_i  \}$ has exact order $h_i$ for all $i$? 
\eenum
\eprob

\section{Sums of finite sets} 

A central problem in additive number theory is to understand the growth of a sumset $hA$ as $h \rightarrow \infty$.

The \emph{$d$-dilation} of the set $A$ is $d\ast A = \{da:a\in A\}$.  
Then $d\ast \N_0 = \{0,d,2d,3d,\ldots\}$ is the set 
of all nonnegative multiples of $d$.   
For infinite sets there is the following result.

\bt[Nash-Nathanson~\cite{nath1985a}] 
Let $A$ be an infinite set of nonnegative integers 
with $d_L(A) > 0$ and  $\gcd(A) = d$.  There is an integer $h_0$ such that $hA \sim d\ast \N_0$
for all $h \geq h_0$.  
\et

For finite sets, a fundamental theorem of additive number theory describes 
the structure of the $h$-fold sum of a finite set $A$ of integers for all sufficiently 
large $h$.  Every finite set of integers 
can be ``normalized'' by a translation and contraction to obtain a set $A$ such that 
$\min(A) = 0$ and $\gcd(A) = 1$.  

\bt[Nathanson~\cite{nath1972-7}]       \label{NSF:theorem:nath1972-7}
Let $A$ be a normalized  finite set of integers with $|A| \geq 2$ 
and $a^* = \max(A)$.  
Let $|A| = k+1$ and 
\[
h^* = (k-1)a^*(a^*-1).
\]
There are integers  $C$ and $D$ and sets $\mcc \subseteq [0,C-2]$ 
and $\mcd \subseteq [0, D-2]$ such that, if $h \geq h^*$, then 
\[
hA = \mcc \cup [C, ha^*-D] \cup (ha^* - \mcd).
\]
\et

Recently,  Granville and Shakan~\cite{gran-shak20} and 
Granville and Walker~\cite{gran-walk21}
obtained the best possible value for $h_0$ and the structure of the 
sets \mcc\ and \mcd\ in this result.  

Theorem~\ref{NSF:theorem:nath1972-7} extends  to linear forms of sums of finite sets of integers.  
A finite sequence of finite sets of nonnegative integers 
can be ``normalized'' by translations and contractions to obtain a 
sequence  $A_1,\ldots, A_r$ such that 
$0 \in A_i$ for all $i \in \{1,\ldots, r\}$
and 
$\gcd\left(\bigcup_{i=1}^r A_i \setminus \{0\}\right) = 1$.

\bt[Han,  Kirfel, and Nathanson~\cite{nath1998-92}]          \label{NSF:theorem:nath1998-92}
Let $A_1,\ldots, A_r$ be a normalized sequence of finite sets of integers with 
$a_i^* = \max(A_i)$ for $i \in \{1,\ldots, r\}$.  There are integers $C$ and $D$ 
and finite sets 
\[
\mcc \subseteq [0,C-2] \qqand \mcd \subseteq [0,D-2]
\]
and there exist integers $h_1^*,\ldots, h_r^*$ such that, 
if $h_i \geq h_i^*$ for all $i \in \{1,\ldots, r\}$, then 
\[
h_1A_1 + \cdots + h_rA_r = \mcc \cup \left[ C, \sum_{i=1}^r h_ia_i^* - D\right] 
\cup \left(  \sum_{i=1}^r h_ia_i^* - \mcd \right). 
\]
\et

\bprob
Determine the best possible values for the integers $h_1^*,\ldots, h_r^*$ 
and the structure of the 
sets \mcc\ and \mcd\ in Theorem~\ref{NSF:theorem:nath1998-92}.  
\eprob

Let $A$ be a set of integers.  
The $h$-fold sumset $hA$ is the set of integers that can be represented 
as the sum of $h$ elements of $A$, that is, the integers $n$ with $r_{h,A}(n) \geq 1$.
For every positive integer $t$, let $(hA)^{(t)}$ be the set of integers 
that have at least $t$ representations  
as  sums of $h$ elements of $A$, that is, the integers $n$ with $r_{h,A}(n) \geq t$.
The structure of the sumsets $(hA)^{(t)}$ has recently been determined. 

\bt[Nathanson~\cite{nath2021-195}]          \label{NSF:theorem:nath2021-195}
Let $A$ be a normalized finite set of integers with $a^* = \max(A)$.  
For every positive integer $t$, let 
\[
h_t^* = (k-1)(ta^*-1)a^* +1.
\]
There are nonnegative integers $C_t$ and $D_t$ and sets $\mcc_t \subseteq [0,C_t-2]$ 
and $\mcd_t \subseteq [0, D_t-2]$ such that, if $h \geq h_t^*$, then 
\[
(hA)^{(t)}= \mcc_t \cup [C_t, ha^*-D_t] \cup (ha^* - \mcd_t).
\]
\et

\bprob
Determine the best possible values for the integers $h_t^*$ 
and the structure of the 
sets $\mcc_t$ and $\mcd_t$ in Theorem~\ref{NSF:theorem:nath2021-195}.  
\eprob

There is a more subtle additive problem.   
Color the elements of the set $A$ with $q$ colors, 
which we call $\{1,\ldots, q\}$.
Let $A_i$ be the subset of $A$ consisting of all elements $a \in A$ that have color $i$. 
The sets $A_1, \ldots, A_q$ are pairwise disjoint with $A = \bigcup_{i=1}^q A_i$.  
The $q$-tuple $\mbA = (A_1,\ldots, A_q)$ is an ordered partition of $A$.  
Let $\N_0^q$ be the set of  $q$-tuples of nonnegative integers.  
For all 
\[
\mbh = (h_1,\ldots, h_q) \in \N_0^q, 
\]  
let
\[
\| \mbh \| = \sum_{i=1}^q h_i = h.
\]
The \emph{chromatic subset}  $\mbh\cdot \mbA$ is the set of all integers in the sumset $hA$ 
that can be represented as the sum of $h$ elements of $A$ 
with exactly $h_i$ elements of color $i$ for all $i \in [1,q]$.  
Thus, $n \in \mbh\cdot \mbA$ if and only if we can write 
\beq                               \label{chromatic:rep-1}
n = \sum_{i=1}^q \sum_{j_i=1}^{h_i} a_{i,j_i}
\eeq
where 
\beq                               \label{chromatic:rep-2}
a_{i,j_i }\in A_i \qquad \text{for all $i \in [1,q]$ and $j_i \in [1,h_i]$.}
\eeq

We refine this problem by allowing elements of $A$ to have more than one color.  
The subset $A_i$ still consists of all elements of $A$ with color $i$, and $A = \bigcup_{i=1}^q A_i$, 
but the sets $A_i$ are not necessarily pairwise disjoint. 
The chromatic sumset $\mbh \cdot \mbA$ consists of all integers $n$ 
that have at least one representation of the form~\eqref{chromatic:rep-1} 
and~\eqref{chromatic:rep-2}. 
Note that 
\[
h_i A_i = \left\{  \sum_{j_i=1}^{h_i} a_{i,j_i}: a_{i,j_i} \in A_i \text{ for all } j_i \in [1,h_i]  \right\} 
\]
and so 
\begin{align*}
\mbh \cdot \mbA & = (h_1,\ldots, h_q)  \cdot (A_1,\ldots, A_q) \\ 
& = h_1A_1 + \cdots + h_qA_q.   
\end{align*}
This is a \emph{homogeneous linear form} in the sets $A_1,\ldots, A_q$.
For every set $B$ of integers, we also have the \emph{inhomogeneous linear form} 
$h_1A_1 + \cdots + h_qA_q + B$.  
Han, Kirfel, and Nathanson (Theorem~\ref{NSF:theorem:nath1998-92}) 
determined the asymptotic structure of homogeneous and inhomogeneous linear forms 
for all $q$-tuples of finite sets of integers.

The \emph{chromatic representation function} $r_{\mbA,\mbh}(n)$ counts the number 
of colored representations of $n$ of the form~\eqref{chromatic:rep-1} and~\eqref{chromatic:rep-2}, 
that is, the number of $q$-tuples 
\[
\left( \mba_1, \ldots, \mba_q \right) 
\in A_1^{h_1}\times \cdots \times \cdots A_q^{h_q} 
\]
where, for all $i \in [1,q]$, the $h_i$-tuple
 $\mba_i  = (a_{i,j_1}, a_{i,j_2}, \ldots,  a_{i,j_{h_i}} )\in A_i^{h_i}$ satisfies 
\[
a_{i,j_1} \leq a_{i,j_2} \leq \cdots \leq a_{i,j_{h_i}}. 
\]
Let 
\[
(\mbh\cdot \mbA)^{(t)} = \{ n \in \mbh\cdot \mbA : r_{\mbA,\mbh}(n) \geq t\}.
\]

Define a partial order on the set $\N_0^q$ as follows.  For vectors 
\[
\mbh_1 = (h_{1,1}, \ldots, h_{q,1}) \in \N_0^q 
\qqand 
\mbh_2 = (h_{1,2}, \ldots, h_{q,2}) \in \N_0^q
\]
let 
\[
\mbh_1 \preceq \mbh_2 \qquad \text{ if $h_{1,i} \leq h_{2,i}$ for all $i \in [1,q]$.}
\]

The following result describes the structure of the chromatic sumset 
$(\mbh\cdot \mbA)^{(t)}$ for all positive integers $t$ and all sufficiently 
large vectors $\mbh  \in \N_0^q$.

\bt[Nathanson~\cite{nath2021-191}]                            \label{NSF:theorem:nath2021-191}
Let $\mbA = (A_1,\ldots, A_q)$ be a normalized $q$-tuple 
of finite sets of integers.  
Let $\max(A_i) = a_i^* \geq 1$ for all $i \in [1,q]$, and
\[
\mba^* = \left( a_1^*,\ldots, a_q^* \right) \in \N_0^q.
\]
For every positive integer $t$, there exist nonnegative integers $c_t$ and $d_t$, 
there exist finite sets of nonnegative integers 
$C_t$ and $D_t$, and there exists a vector $\mbh_t = (h_{t,1},\ldots, h_{t,q} ) \in \N_0^q$ such that, 
if  $\mbh = (h_{1},\ldots, h_q) \in \N_0^q$ and $\mbh \succeq \mbh_t$, then 
\[
\left(\mbh \cdot \mbA\right)^{(t)}  = C_t \cup \left[ c_t, \mbh \cdot \mba^* - d_t  \right] 
\cup \left( \mbh \cdot \mba^* - D_t \right). 
\] 
\et

\section{Lattice points and abelian semigroups}

Let $\Z^n$ denote the group of lattice points in $\R^n$. 
The higher-dimensional additive problem is the study of sums of finite sets of lattice points.  
In discrete geometry, a lattice polytope is the convex hull of a finite set of lattice points.  
Let $A$ be  a finite subset of $\Z^n$ that is not contained in a hyperplane of $\R^n$, 
and let $K$ be the convex hull of $A$.
For every positive integer $h$, the dilation $h\ast K$ is equal to the $h$-fold sumset $hK$.
Ehrhart ~\cite{ehrh67a,ehrh67b,ehrh68} proved that there is a polynomial $e_K(z)$ of degree $n$ 
such the number of lattice points in 
the sumset $hK$, that is, $\card(hK \cap \Z^n)$, 
is exactly $e_K(h)$ for all sufficiently large $h$.  
We have $A \subseteq K$ and so $hA \subseteq hK \cap \Z^n$.  It follows that 
\[
\card(hA) \leq \card(hK \cap \Z^n) = e_K(h).
\]
Khovanskii~\cite{khov92,khov95} 
proved that there is a polynomial $f_A(z)$ of degree $n$ 
such that 
\[
\card(hA \cap \Z^n)  = f_A(h)
\]
 for all sufficiently large $h$ and that the sumset $hA$ contains every lattice 
 point of $hK$ that is not ``close'' to the boundary of the polytope $hK$.  
 
 \bprob
 Describe the structure of the set of lattice points in the ``boundary layer'' 
 of points in $(hK \cap \Z^n) \setminus hA$. 
 \eprob

Now let $X$ be any additive abelian semigroup.  
Khovanskii  
proved that if $A$ and $B$ are finite subsets of $X$, 
then there is a polynomial $g_{A,B}(z)$ 
such that 
\[
\card(hA + B)   = g_{A,B}(h) 
\]
 for all sufficiently large $h$.

Nathanson extended this result to linear forms of finite sequences of finite subsets of $X$.

\bt[Nathanson~\cite{nath2000-98}]
Let $A_1,\ldots, A_{\ell}$, and $B$ be finite subsets 
of an additive abelian semigroup $X$.  
There exists a polynomial $p_{A_i,B}(z_1,\ldots, z_{\ell})$ 
and an integer $h^* = h^*(A_1,\ldots, A_{\ell},B)$ 
such that, if $h_1,\ldots, h_{\ell}$ are integers with $\min(h_i) \geq h^*$, 
then
\[
\card\left( h_1A_1+\cdots + h_{\ell}A_{\ell} + B ) \cap \Z^n \right)  = p_{A_i,B}(h_1,\ldots, h_{\ell}) .
\]
\et

The proof uses the Hilbert polynomial for finitely generated graded rings.
Nathanson and Ruzsa~\cite{nath2002-101}  
gave an elementary proof of the special case  $B=\{0\}$. 

 \section{Approximate groups} 
Let $A$ be a nonempty subset of an additive abelian group $G$, 
and let $r$ and $\ell$ be positive integers.  
The set $A$ is an \emph{$(r,\ell)$-approximate group} 
if there exists a set $X \subseteq G$ such that 
\beq     \label{AAG-def1}
|X| \leq \ell
\eeq
and
\beq     \label{AAG-def2}
rA\subseteq X+A.
\eeq
The idea of an approximate group evolved from Freiman's inverse theorem 
in additive number theory (Freiman~\cite{frei64,frei64t,frei66,frei73}, 
Nathanson~\cite{nath1996bb}).
This definition of approximate group is less restrictive  
than the original definition, which appears in Tao~\cite{tao07} 
and which has been extensively investigated (for example, by 
Breuillard, Green, and Tao~\cite{breu13,breu15,breu-gree-tao12,gree12,gree14} 
and Pyber and Szabo~\cite{pybe-szab14}).  
Helfgott~\cite{helf15} is a recent survey of related problems on growth in groups.

The set $A$ is an \emph{asymptotic $(r,\ell)$-approximate group} 
if, for all sufficiently large $h$,  
the $h$-fold sumset of $A$ is an $(r,\ell)$-approximate group.
This means that there exists an integer $h_0(A)$ such that, for each integer
 $h \geq h_0(A)$, 
there is a set $X_h \subseteq G$ such that 
\beq     \label{AAG-def1-asymp}
|X_h| \leq \ell
\eeq
and
\beq     \label{AAG-def2-asymp}
rhA \subseteq X_h + hA.
\eeq

\bt[Nathanson~\cite{nath2018-179}]           \label{nsf:theorem:nath2018-179}
Every nonempty finite subset of an abelian group is an asymptotic approximate group. 
\et

\bprob
Theorem~\ref{nsf:theorem:nath2018-179} is a very general result.  It would 
be of interest to obtain explicit estimates for $r$, $h$, and $\ell$ for 
finite sets of lattice points.  
\eprob

\section{Sidon sets and $B_h$-sets}

The Erd\H os-Tur\' an conjecture asserts the unboundedness of the representation 
function of an asymptotic basis.  In the opposite direction, we consider sets 
with sumsets whose representation functions are bounded.  The extreme case 
 is a \emph{Sidon set}, which is a finite or infinite subset $A$ of the 
 nonnegative integers (or of any additive abelian semigroup) 
such that $r_{2,A}(n) \leq 1$ for all $n \in \N_0$.  
More generally,  a  \emph{$B_h$-set} is a finite or infinite set $A$ of 
the nonnegative integers 
such that $r_{h,A}(n) \leq 1$ for all $n \in \N_0$.  A Sidon set is a  $B_2$-set.
A \emph{$B_h[g]$-set} is a finite or infinite set $A$ of nonnegative integers 
such that $r_{h,A}(n) \leq g$ for all $n \in \N_0$.  
There is a vast literature on these sets (O'Bryant~\cite{obry04}) and many open problems.  
For example, what is the cardinality of the largest Sidon set contained in the interval $[1,n]$? 
The current best result is due to Balogh,  F\"{u}redi, and Roy~\cite {balo-fure-roy23}.

Of recent interest are perturbations of Sidon sets and $B_h$-sets.
Let \mbF\ be a field with a nontrivial absolute value $| \ |$. 
Let $I$ be a nonempty set, and let $\{a_i : i\in I \}$ and $\{b_i : i\in I \}$ 
be sets of elements of the field \mbF.  
Let $\{\varepsilon_i : i\in I \}$ be a set of positive real numbers.  
The set  $\{b_i : i\in I \}$  is an  
\emph{$\varepsilon$-perturbation}\index{perturbation}  
of the set $\{a_i : i\in I \}$ if 
$|b_i - a_i| < \varepsilon_i$
for all $i \in I$.

\bt[Nathanson~\cite{nath2021-196}]         \label{perturb:theorem:InfiniteSidon} 
Let \mbF\ be a field of characteristic 0 with a nontrivial absolute value.   
Let $\{\varepsilon_i : i \in \N\}$ be a set of positive real numbers.  
For every subset $\{a_i: i \in \N \}$ of \mbF, there is a $B_h$-set $\{b_i: i \in \N \}$ in \mbF\  
that is an $\epsilon$-perturbation of $A$.
\et

\bprob
Let $A$ be an uncountably infinite set.  Is there an $\varepsilon$-perturbation of $A$ that is an $h$-Sidon set? 
\eprob

The \emph{greedy $B_h$-set} is the $B_h$-set constructed inductively 
by letting $a_0(h) = 0$, $a_1(h) = 1$, and, given the $B_h$-set 
$\{a_0(h) , a_1(h) ,\ldots, a_k(h) \}$,
choosing $a_{k+1}(h)  > a_k(h) $ as the smallest integer such that 
$\{a_0(h) , a_1(h) ,\ldots, a_k(h) , a_{k+1}(h) \}$ is a $B_h$-set.  
We know very little  about the structure of these much investigated, but still mysterious, 
sets constructed with the greedy algorithm. 
We have $a_0(h) =0$, $a_1(h) = 1$, and $a_2(h) = h+1$.  
Recently, Nathanson~\cite{nath2023-220} proved that 
\[
a_3(h) = h^2+h+1
\]
for all $h \geq 1$.
Nathanson and O'Bryant~\cite{nath2023-221}  
proved that 
\[
a_4(h) = \begin{cases}
\left( h+1\right)^3/2 & \text{if $h$ is odd}\\
\left(  h^3 + 2h^2 +3h+2\right)/2& \text{if $h$ is even.}
\end{cases}
\]
O'Bryant~\cite{obry24} has determined upper and lower bounds for $a_5(h)$.

\bprob
Compute $a_5(h)$ and determine if $a_k(h)$ 
is a quasi-polynomial 
for all $h \geq 5$.
\eprob

For other recent work on Sidon sets, see Nathanson~\cite{nath2021-196, nath2022-203, nath2022-200}.

\section{Sidon sets for linear forms}

Let \mbF\ be a field and let $h$ be a positive integer.  
We consider linear forms  
\beq                                                      \label{LinearPerturb:phi}
\varphi(x_1,\ldots, x_h) =  c_1 x_1 + \cdots + c_h x_h  
\eeq
where $c_i \in \mbF$ for all $i \in \{1,\ldots, h\}  $.  
Let $A$ be  a nonempty subset of  a vector space $V$ over the field \mbF.  
The \emph{$\varphi$-image of $A$} is the  sum of dilates 
\begin{align*}
\varphi(A) & = \left\{ \varphi(a_{1},\ldots, a_{h}): (a_1,\ldots, a_h) \in A^h \right\} \\ 
& = \left\{ c_1 a_{1} + \cdots +  c_h a_{h}: (a_1,\ldots, a_h) \in A^h \right\} \\
& = c_1\ast A + \cdots + c_h \ast A. 
\end{align*}
A nonempty subset $A$  of $V$ is a 
\emph{Sidon set for the linear form $\varphi$} or, simply, a 
\emph{$\varphi$-Sidon set} 
if it satisfies the following property:  
For all $h$-tuples $(a_1,\ldots, a_h) \in A^h$ and $ (a'_1,\ldots, a'_h) \in A^h$,   
if 
\[
\varphi(a_1,\ldots, a_h) = \varphi(a'_1,\ldots, a'_h) 
\]
then $(a_1,\ldots, a_h) = (a'_1,\ldots, a'_h)$, 
that is, $a_i = a'_i$ for all $i \in \{1,\ldots, h\} $.

There is a simple obstruction to the existence of $\varphi$-Sidon sets 
 with more than one element.  For every nonempty subset $I$ of $\{1,\ldots, h\}$, 
 define the \emph{subset sum} 
 \beq                                                            \label{LinearPerturb:sI}
 s(I) = \sum_{i \in I} c_i. 
\eeq
Let  $s({\emptyset}) = 0$.  
Suppose there exist disjoint subsets $I_1$ and $I_2$ of $\{1,\ldots, h\}$ 
with $I_1$ and $I_2$ not both empty such that 
\beq                                                      \label{LinearPerturb:obstruction}
s({I_1}) = \sum_{i \in I_1} c_{i} = \sum_{i \in I_2} c_{i} = s({I_2}).   
\eeq
Let $I_3 = \{1,\ldots, h\}  \setminus (I_1 \cup I_2)$.  
Let $A$ be a subset of $V$ with  $|A| \geq 2$.  Choose vectors $u,v,w  \in A$ with $u \neq v$, 
and define
\[
a_i = \begin{cases}
u & \text{ if $i \in I_1$}\\
v & \text{ if $i \in I_2$} \\
w & \text{ if $i \in I_3$}
\end{cases}
\qqand
a'_i = \begin{cases}
v & \text{ if $i \in I_1$ } \\
u & \text{ if $i \in I_2 $}  \\
w & \text{ if $i \in I_3$.}
\end{cases}
\]
We have 
\[
(a_1,\ldots, a_h) \neq (a'_1,\ldots, a'_h)
\]
because $I_1 \cup I_2 \neq \emptyset$ and $a_i \neq a'_i$ for all $i \in I_1 \cup I_2$.  

The sets $I_1$, $I_2$, $I_3$ are pairwise disjoint.   
Condition~\eqref{LinearPerturb:obstruction} implies 
\begin{align*}
\varphi(a_1,\ldots, a_h) 
& = \left(  \sum_{i\in I_1} c_i \right) u  + \left( \sum_{i\in I_2} c_i  \right) v +  \left( \sum_{i\in I_3} c_i \right) w \\
& = \left( \sum_{i\in I_2} c_i \right) u  +  \left( \sum_{i\in I_1} c_i  \right) v+  \left( \sum_{i\in I_3} c_i \right) w \\
& = \varphi(a'_1,\ldots, a'_h)  
\end{align*} 
and so $A$ is not a $\varphi$-Sidon set.  
We say that the linear form~\eqref{LinearPerturb:phi} has \emph{property $N$}  
if there do \emph{not} exist disjoint subsets $I_1$ and $I_2$ of  $\{1,\ldots, h\}$ 
that satisfy condition~\eqref{LinearPerturb:obstruction} 
with $I_1$ and $I_2$ not both empty.  

Sidon sets for linear forms were introduced by Nathanson, 
who proved the following theorems.

\bt[Nathanson~\cite{nath2022-202}]                                          \label{LinearPerturb:theorem:exist-V}
Let \mbF\ be a field, let $V$ be an infinite vector space over the field \mbF, 
and let $X$ be an infinite subset of $V$.   
Let $\varphi(x_1,\ldots, x_h) = \sum_{i=1}^h c_ix_i$ be a linear form 
with nonzero coefficients $c_i \in \mbF$.   
The following are equivalent:
\benum
\item[(i)] 
The set $X$ contains an infinite $\varphi$-Sidon set $A$.
\item[(ii)]
The set $X$ contains a $\varphi$-Sidon set $A$ with $|A| \geq 2$. 
\item[(iii)]
The linear form $\varphi$ has property $N$.  
\eenum
\et

\bt[Nathanson~\cite{nath2022-202}]                                \label{LinearPerturb:theorem:asymptotic}
Let \mbF\ be a field  with a nontrivial absolute value  
and let $V$ be a vector space over \mbF\ that has a norm with respect to the absolute value 
on \mbF.  Let $\varepsilon = \{\varepsilon_k: k = 1,2,3,\ldots \}$ be a set
of positive real numbers.
Let $\varphi$ be a linear form with coefficients in \mbF\ that has property $N$.  
For every set $B = \{ b_k: k = 1,2,3,\ldots \}$ of vectors in $V$, 
there is a $\varphi$-Sidon set $A= \{a_k: k = 1,2,3,\ldots \}$ of vectors in $V$ 
such that 
 \beq                                                                                   \label{LinearPerturb:inequality}
\| a_k - b_k \| < \varepsilon_k 
\eeq
 for all $k = 1,2,3,\ldots$. 
\et

\bt[Nathanson~\cite{nath2022-202}]                                                       \label{LinearPerturb:theorem:limit}
Let \mbF\ be a field  with a nontrivial absolute value, 
and let $\varphi$ be a linear form with coefficients in \mbF\ that has property $N$.  
Let $V$ be a vector space over \mbF\ that has a norm with respect to absolute value on \mbF. 
For every set $B = \{b_k: k=1,2,3,\ldots \}$ of vectors in $V$, there exists a $\varphi$-Sidon set
 $A = \{a_k: k=1,2,3,\ldots \}$ in $V$ such that 
 \[
\lim_{k \rightarrow \infty} \| a_k  - b_k \| = 0.
 \]
\et

Let $\PP = \{2,3,5,\ldots \}$ be the set of prime numbers.  
For every prime number $p$, let $| \ |_p$ be the usual $p$-adic absolute value with $|p|_p = 1/p$.

\bt[Nathanson~\cite{nath2022-202}]        \label{LinearPerturb:theorem:Positive}
Let $\varphi$ be a linear form with rational coefficients that satisfies property $N$.  
Let $\{\varepsilon_k : k=1,2,3,\ldots \}$ be a sequence of positive real numbers 
and let $\{p_k:k=1,2,3,\ldots\}$ be a sequence of prime numbers.    
For every sequence of integers $B = \{b_k:k=1,2,3,\ldots\}$, 
there exists a strictly increasing sequence of positive integers $A =\{ a_k:k=1,2,3,\ldots\}$ such that 
$A$ is a $\varphi$-Sidon set and 
\[
 | a_k - b_k |_{p_j} < \varepsilon_k 
\]
for all $k \in \N$ and $j \in \{1,\ldots, k\}$.  
\et

\bprob
Let $\varphi = \sum_{i=1}^h c_i x_i$ be a linear form with integer coefficients.  
Let $\PP$ be the set of prime numbers and let $A = \{\log p: p \in \PP\}$.
Consider the $h$-tuple $(p_1,\ldots, p_h) \in \PP^h$ of not necessarily distinct prime numbers, 
and let $\PP_0 = \{p \in \PP : p=p_i \text{ for some } i \in \{1,\ldots, h\} \}$.   
For each $p \in \PP_0$, let 
\[
I_p = \{i\in \{1,\ldots, h\}: p_i = p \} \qqand s(I_p) = \sum_{i\in I_p} c_i.
\]
We have 
\[
\varphi(\log p_1,\ldots, \log p_h) = \sum_{i=1}^h c_i \log p_i = \sum_{p\in \PP_0} s(I_p) \log p 
= \log \prod_{p\in \PP_0} p^{S(I_p)}. 
\]
If the linear form $\varphi$ satisfies property $N$, then, by the fundamental theorem of arithmetic, 
the set $A = \{\log p: p \in \PP\}$ is a $\varphi$-Sidon set.  

For the linear form $\psi = x_1 + \cdots + x_h$, Ruzsa~\cite{ruzs98a} used the set $A$ to construct 
large classical Sidon sets of positive integers .  
Are such constructions also possible for $\varphi$-Sidon sets of positive integers? 
\eprob

\bprob 
Let $A = \{ a_k : k=1,2,3,\ldots \}$ and $B = \{ b_k: k = 1,2,3,\ldots\}$ 
be sequences of integers.  
The set $A$ is a \emph{polynomial perturbation} of  $B$ if 
\[
|a_k - b_k| < k^r 
\] 
for some $r > 0$ and  all $k \geq k_0$.  
The set $A$ is a \emph{bounded perturbation} of  $B$ if 
\[
|a_k - b_k| < m_0
\]
for some $r > 0$ and  all $k \geq k_0$.  

Let $\varphi$ be a linear form with integer coefficients that satisfies condition $N$.  
Let $B$ be a set of integers.  Does there exist a $\varphi$-Sidon set of integers 
that is a polynomial perturbation of $B$?  

Does there exist a $\varphi$-Sidon set of integers 
that is a bounded perturbation of $B$?

\eprob

\bprob 
Let $\varphi$ be a linear form with integer coefficients that satisfies condition $N$.  
For every positive integer $n$, determine the cardinality of the largest $\varphi$-Sidon 
subset of $\{1,2,\ldots, n\}$.

\eprob

\bprob 
There exists $c > 0$ such that, for every positive integer $n$, 
 there is a classical Sidon set $A \subseteq \{1,\ldots, n\}$ 
with $A(n) \geq c  \sqrt{n}  $.  
However, there is no infinite classical Sidon set $A$ of positive integers 
such that $A(n) \geq c \sqrt{n} $ for some $c > 0$ and all $n \geq n_0$.  
Indeed, Erd\H os (in St{\"o}hr~\cite{stoh55}) proved that 
every infinite classical Sidon set satisfies 
\[
\liminf_{n\rightarrow \infty} A(n)\sqrt{ \frac{\log n}{n} } \ll 1.
\]
Are there analogous lower bounds for infinite $\varphi$-Sidon sets of positive integers 
associated with binary linear forms $\varphi = c_1x_1 + c_2x_2$ or with linear forms 
$\varphi = \sum_{i=1}^h c_i x_i$ for $h \geq 3$? 

\eprob

\bprob 
Consider  sets of integers.  
One might expect that the elements of a set $A$ of integers that is ``sufficiently random'' or 
``in general position'' will be a  classical Sidon set, that is, will not contain a nontrivial solution 
of the equation $x_1+x_2 = x_3+x_4$.  Equivalently, the set $A$ will be one-to-one (up to transposition) on 
the function $f(x_1,x_2) = x_1 + x_2$.  There is nothing special about the function $x_1 + x_2$.  
One could ask if $A$ is one-to-one  (up to permutation) on some symmetric function, 
or one-to-one on a function that is not symmetric.  
The  functions considered in this paper are linear forms in $h$ variables.  

Conversely, given the set $A$ of integers, we can ask what are the functions 
(in some particular set \mcf\ of functions) with respect to which the set $A$ is one-to-one.  
This inverse problem is considered in Nathanson~\cite{nath2022-202}.

\eprob

\def\cprime{$'$} \def\cprime{$'$} \def\cprime{$'$}
\providecommand{\bysame}{\leavevmode\hbox to3em{\hrulefill}\thinspace}
\providecommand{\MR}{\relax\ifhmode\unskip\space\fi MR }
\providecommand{\MRhref}[2]{%
  \href{http://www.ams.org/mathscinet-getitem?mr=#1}{#2}
}
\providecommand{\href}[2]{#2}

\end{document}